\documentclass[a4paper,12pt,twoside,notitlepage]{article}
\usepackage{a4,amsmath,amssymb}
\usepackage[noblocks]{authblk}
\usepackage{graphicx}
\def\q0{\theta}
\def\q{\vartheta}

\def\e0{\epsilon}

\def\f0{\phi}
\def\f{\varphi}

\def\rd{\mathrm d}

\def\R{{\mathbb R}}
\def\M{{\mathbf M}}

\def\pab{\partial}

\def\e{\varepsilon}

\newtheorem{thm}{Theorem}[section]
\newtheorem{cor}[thm]{Corollary}

\newtheorem{prop}[thm]{Proposition}
\newtheorem{defn}[thm]{Definition}
\newtheorem{rem}[thm]{Remark}

\begin{document}
\title{\sc The discrete centroaffine indefinite surface}
\author{Yun Yang, Yanhua Yu\thanks{Corresponding author. \\ {\it \hspace*{1em} Email addresses}: freeuse\_st@126.com (Yun Yang), yyh\_start@126.com(Yanhua Yu).}
\\{\it\small Department of Mathematics, Northeastern University,}\\{\it\small Shenyang, Liaoning, P. R. China,  110004}}
\markboth{}{}
\date{}
\maketitle
\numberwithin{equation}{section}
\begin{abstract}
In this paper we build the structure equations and the integrable systems for a discrete centroaffine indefinite surface in $\R^3$. At the same time, some centroaffine invariants are obtained according to the structure equations. Using these centroaffine invariants, we study the Laplacian operator and the convexity of a discrete centroaffine indefinite surface. Furthermore, some interest examples are provided.
\end{abstract}
\section{Introduction.}
Discrete differential geometry studies discrete equivalents of the geometric notions and methods of classical differential
geometry, such as notions of curvature and integrability for polyhedral surfaces.  In this connection, discrete surfaces have been studied one after another with
strong ties to mathematics physics and great potential for computer analysis, architecture, numerics. Progress in this field is to a
large extent stimulated by its relevance for computer graphics and mathematical physics\cite{Bobenko-2,Bobenko-4}.
\par Recently, the expansion of computer graphics and applications in mathematical physics have given a great impulse to the issue of giving discrete equivalents of affine differential geometric objects\cite{Bobenko-1,Craizer-1,Craizer-2}. In \cite{Bobenko-3} a consistent definition of discrete affine spheres is proposed, both for definite and indefinite metrics, and in \cite{Matsuura} a similar construction is done in the context of improper affine spheres.
\par Following the ideas of Klein, presented in his famous lecture at Erlangen, several geometers in the early 20th century proposed the study
of curves and surfaces with respect to different transformation groups. In geometry, an affine transformation, affine map or an affinity is a function between affine spaces which preserves points, straight lines and planes. Also, sets of parallel lines remain parallel after an affine transformation. An affine transformation does not necessarily preserve angles between lines or distances between points, though it does preserve ratios of distances between points lying on a straight line. Examples of affine transformations include translation, scaling, homothety, similarity transformation, reflection, rotation, shear mapping, and compositions of them in any combination and sequence.
\par A centroaffine transformation is nothing but a general linear transformation $\R^n\ni x\mapsto Ax\in\R^n$, where $A\in GL(n,\R)$.
 In 1907 Tzitz$\mathrm{\acute{e}}$ica found that for a surface in Euclidean 3-space the property that the ratio of the Gauss curvature to the fourth power of
the distance of the tangent plane from the origin is constant is invariant under a centroaffine transformation. The surfaces with this property
turn out to be what are now called Tzitz$\mathrm{\acute{e}}$ica surfaces, or proper affine spheres with center at the origin. In centroaffine differential geometry, the theory of hypersurfaces has a long history. The notion of centroaffine minimal
hypersurfaces was introduced by Wang \cite{Wang} as extremals for the area integral of the centroaffine metric. See also \cite{Y-Y-L,Yu-Y-L} for the classification results about centroaffine translation surfaces and centroaffine ruled surfaces in $\R^3$.
\par Smooth geometric objects and their transformations should belong to the same geometry. In particular discretizations should be invariant with respect to the same transformation group as the smooth objects are(projective, affine, m\"obius etc). This paper is concerned with some invariant properties of the discrete centroaffine indefinite surface, which is organized as follows: Basic concepts of classical centroaffine differential geometry are presented in Section 2.
In Section 3 we define the discrete centroaffine indefinite surface, and then obtain the structure equations, compatibility conditions and some centroaffine invariants.
In section 4, the Laplacian operator is defined as the gradient of the Dirichlet energy for a discrete centroaffine indefinite surface. The discrete centroaffine indefinite surface with constant coefficients is considered in section 5. Section 6 deals with the convexity of a discrete centroaffine indefinite surface.
\section{Centroaffine hypersurfaces.}
Prior to the introduction of a discrete centroaffine indifinite surfaces theory, in this section, we recall some fundamental notions for centroaffine
hypersurfaces in $\R^{n+1}$. For details we refer to
\cite{L-S-W}, \cite{Liu-W-1}, \cite{S-S-V} or \cite{Wang}. Let $x: \M \rightarrow
\R^{n+1}$ be a hypersurface immersion and $[\dots]$ the standard
determinant in $\R^{n+1}$. $x$ is said to be a centroaffine
hypersurface if the position vector of $x$, denoted also by $x$, is
always transversal to the tangent space $ x_*({\mathbf{T}\M})$ at
each point of $\M$ in $\R^{n+1}$. We define a symmetric bilinear
form $G$ on ${\mathbf T}\M$ by
\begin{equation}\label{metric}
G=-\sum_{i,j=1}^n{\frac{[e_1(x),e_2(x),\dots,e_n(x),e_ie_j(x)]}{[e_1(x),e_2(x),\dots,e_n(x),x]}}\
\theta ^{i}\otimes \theta ^j,
\end{equation}
where $\{ e_1,e_2,\dots,e_n \}$ is a local basis of $\mathbf T\M$ with
the dual basis $\{\theta^1,\theta^2,\dots,\theta^n\}$. Note that $G$ is
globally defined. A centroaffine hypersurface $x$ is said to be non-degenerate if $G$ is non-degenerate.
We call $G$ the centroaffine metric of $x$.
We say that a hypersurface is definite (or indefinite) if $G$ is definite (or indefinite).
\begin{rem}
Geometrically, a hypersurface $x$ with positive (resp. negative) definite centroaffine metric $G$ is the
locally strongly convex hypersurface in
$\R^{n+1}$ and such hypersurface is called hyperbolic type (respectively,  elliptic type) in \cite{L-L-S}.
In particular, $G$ is definite if $x$ is locally
strongly convex in $\R^{n+1}$.
\end{rem}
\par
Let $x: \M \rightarrow \R^{n+1}$ be a non-degenerate centroaffine
surface. Then $x$ induces a centroaffinely invariant metric $G$ and
a so-called induced connection $\nabla$. The difference of the
Levi-Civita connection $\widehat{\nabla}$ of $G$ and the induced
connection $\nabla$ is a $(1,2)-$tensor $C$ on $M$ with the property
that its associate cubic form $\widehat{C}$, defined by
\begin{equation}
\widehat{C}(u,v,w)=G(C(u,v),w), u,v,w\in TM,
\end{equation}
which is totally symmetric. The so-called Tchebychev form is defined by
\begin{equation}\label{Tch-form}
\widehat{T}=\frac{1}{n}\mathrm{trace}_{G}(\widehat{C}).
\end{equation}
Let $T$ be the Tchebychev vector field on $M$ defined by the equation
\begin{equation}\label{Tch-vector}
G(T,v)=\widehat{T}(v), v\in TM.
\end{equation}

It is proved by Wang in \cite{Wang} that a centroaffine surface $x:\M^n \rightarrow \R^{n+1}$ is called centroaffine minimal if $\mathrm{trace}_{G}(\widehat{\nabla}T)=0$,
and the centroaffine mean curvature is defined by
\begin{equation}\label{H-def}
H=\frac{1}{n}\mathrm{trace}_{G}(\widehat{\nabla}T).
\end{equation}
The Gauss equation of $x$ can be written as(in the following, we use
the Einstein summation convention and the range of indices is $1\leq
i,j,k,\dots\leq n$)
\begin{equation}\label{gauss-eqn}
    \frac{\pab^2x}{\pab x^i\pab x^j}=\Gamma_{ij}^ke_k(x)-G_{ij}x.
\end{equation}
\par Then the Riemannian curvature tensor is given by
\begin{equation}
    \widehat{R}^{l}_{ijk}=\frac{\pab\widehat{\Gamma}^{l}_{ij}}{\pab u^k}-\frac{\pab\widehat{\Gamma}^{l}_{ik}}{\pab u^j}+\widehat{\Gamma}^{p}_{ij}\widehat{\Gamma}^{l}_{pk}
    -\widehat{\Gamma}^{p}_{ik}\widehat{\Gamma}^{l}_{pj} ,
\end{equation}
and
\begin{equation}
\widehat{R}_{mijk}=G_{ml}\widehat{R}^{l}_{ijk},
\end{equation}
where $\widehat{\Gamma}^{k}_{ij}$ is the Levi-Civita connection of $G$.
\par If $n=2$, the Gauss curvature of $x$ is defined by
\begin{equation}\label{gauss-cuv-def}
\kappa=\frac{-\widehat{R}_{1212}}{\det(G_{ij})}.
\end{equation}

\par Let $x:M\rightarrow \mathbb{R}^3$ be an centroaffine indefinite surface. We introduce local asymptotic coordinates $(u, v)$ of $G$ such that
\begin{equation}
  G=h(\rd u\otimes\rd v+\rd v\otimes\rd u)
\end{equation}
for some local function $h>0$. Using appropriate functions $\lambda,\mu,\varphi,\psi$ we define $1$-forms
\begin{equation}
  \Lambda:=\lambda\rd u:=\frac{[x_{uu},x_u,x]}{[x_u,x_v,x]}\rd u,\quad \mathcal{M}:=\mu\rd v=\frac{[x_{uv},x_v,x]}{[x_u,x_v,x]}\rd v
\end{equation}
and cubic forms
\begin{equation}
  \Phi:=\varphi\rd u^3:=h\frac{[x_u,x_{uu},x]}{[x_u,x_v,x]}\rd u^3,\quad \Psi:=\psi\rd v^3=-h\frac{[x_{u},x_{vv},x]}{[x_u,x_v,x]}\rd v^3.
\end{equation}
Then we have the following structure equations
\begin{eqnarray}
  x_{uu} &=& (\frac{h_u}{h}+\lambda)x_u+\frac{\varphi}{h}x_v, \\
  x_{uv} &=& \mu x_u+\lambda x_v-hx, \\
  x_{vv} &=& \frac{\psi}{h}x_u+(\frac{h_v}{h}+\mu)x_v
\end{eqnarray}
and the integrability conditions
\begin{eqnarray}
  &(\ln h)_{uv}+\frac{\varphi\psi}{h^2}-\lambda\mu+h = 0 \\
  &\varphi_v=h(\lambda_u-\frac{h_u}{h}\lambda),\quad \psi_u=h(\mu_v-\frac{h_v}{h}\mu), \\
  &\lambda_v=\mu_u.
\end{eqnarray}
We will have
\begin{eqnarray}
 &C^1_{11}=\frac{h_u}{2}+\lambda,\quad C^2_{11}=\frac{\varphi}{h},\\
 &C^1_{12}=\mu,\quad C^2_{12}=\lambda,\\
 &C^1_{22}=\frac{\psi}{h},\quad C^2_{22}=\frac{h_v}{h}+\mu,\\
 &T^1=\frac{\mu}{h},\quad T^2=\frac{\lambda}{h}.
\end{eqnarray}

Let $\ll,\gg$ be the inner product of the forms on $M$ induced by the centroaffine metric $G$. Then, by the definition of the centroaffine metric, the Tchebychev form, and the forms $\Lambda, \mathcal{M}, \Phi$ and $\Psi$, we have
\begin{eqnarray}
  &\ll \Lambda, \mathcal{M}\gg=\frac{\lambda\mu}{h} =\frac{1}{2}||T||^2,\\
  &\ll \Phi, \Psi\gg=\frac{\varphi\psi}{h^3}=J-\frac{3}{2}||T||^2=\tilde{J}, \\
  & \kappa=-\frac{(\ln h)_{uv}}{h},\\
  & \kappa-1=\ll \Phi, \Psi\gg-\ll \Lambda, \mathcal{M}\gg,
\end{eqnarray}
where, $J$ is the Pick invariant and $\kappa$ the Gauss curvature of $G$. There are the
following propositions:
\begin{prop}(i) $\Lambda=\mathcal{M}=0$ if and only if x is a proper equiaffine sphere centered at
origin $O\in\mathbb{R}^3$; (ii) $\Phi=\Psi=0$ if and only if x is a quadric (\cite{Liu-W-1}).
\end{prop}
\begin{prop}
(i) An centroaffine indefinite surface $x$ in $\R^3$ is a centroaffine extremal
(minimal) surface if and only if $\Lambda$ and $\mathcal{M}$ satisfy $ \lambda_u= \mu_v = 0$; (ii) an indefinite centroaffine
surface $x$ in $\R^3$ is a centroaffine Tchebychev surface if and only if $\Phi$ and $\Psi$	 satisfy $\varphi_v=\psi_u=0 $(\cite{Liu-W-1}).
\end{prop}

\section{Discrete centroaffine indefinite immersions.}
Here, we define discrete analogues of centroaffine immersions in a purely geometric manner. These constitute particular `discrete surface' which are maps
\begin{equation}\label{dis-map}
  \vec{r}:\mathbb{Z}^2\rightarrow \mathbb{R}^3,\qquad (n_1,n_2)\mapsto \vec{r}(n_1,n_2).
\end{equation}
In the following, we suppress the arguments of functions of $n_1$ and $n_2$ and denote increments of the discrete variables by subscripts, for example,
$$\vec{r}=\vec{r}(n_1,n_2),\vec{r}_1=\vec{r}(n_1+1,n_2),\vec{r}_2=\vec{r}(n_1,n_2+1).$$
Moreover, decrements are indicated by overbars, that is,
$$\vec{r}_{\bar{1}}=\vec{r}(n_1-1,n_2),\vec{r}_{\bar{2}}=\vec{r}(n_1,n_2-1).$$
The following notation for difference operators is adopted:
$$\Delta_i\vec{r}=\vec{r}_i-\vec{r},\quad \Delta_{12}=\vec{r}_{12}-\vec{r}_1-\vec{r}_2+\vec{r}.$$
Now we will give a definition for the discrete centroaffine indefinite surface. Especially, in the following, the point $\vec{r}(n_1,n_2)$ indicates the terminal point of the vector $\vec{r}(n_1,n_2)$ with its starting point at the origin $O$.
\begin{defn}\label{dics}(Discrete centroaffine indefinite surface) A two-dimensional lattice (net) in three-dimensional affine space
\begin{equation}
 \vec{r}:\mathbb{Z}^2\rightarrow \mathbb{R}^3
\end{equation}
is called a discrete centroaffine indefinite surface if it has the following properties:
\begin{enumerate}
  \item[(a)] Any point $\vec{r}(n_1,n_2)$ and its neighbours $\vec{r}_{\bar{1}},\vec{r}_{\bar{2}},\vec{r}_1,\vec{r}_2$ lie on a plane $\pi$.
  \item[(b)] The origin $O$ is not in the plane $\pi$.
  \item[(c)] The three points $\vec{r}(n_1,n_2), \vec{r}(n_1+\epsilon_1,n_2), \vec{r}(n_1,n_2+\epsilon_2)$ are nonlinear, where $\epsilon_1,\epsilon_2\in\{1,-1\}$.
\end{enumerate}
\end{defn}
In analytical terms, Condition (a) can be translated into
\begin{equation}
\det[\vec{r}_1-\vec{r},\vec{r}_2-\vec{r},\vec{r}-\vec{r}_{\bar{1}}]=0,\quad \det[\vec{r}_1-\vec{r},\vec{r}_2-\vec{r},\vec{r}-\vec{r}_{\bar{2}}]=0.
\end{equation}
Condition (b) and Condition (c) imply
\begin{equation}\label{ind-vec}
  \det[r_1,r_2,r]\neq0,\quad \det[r_{\bar{1}},r_{\bar{2}},r]\neq0,\quad \det[r_1,r_{\bar{2}},r]\neq0,\quad \det[r_{\bar{1}},r_2,r]\neq0.
\end{equation}
\par So that the position vector of the discrete surface considered here obeys the discrete `Gauss equation'
\begin{eqnarray}
  \vec{r}_{11}-\vec{r}_1 &=& \alpha(\vec{r}_1-\vec{r})+\beta(\vec{r}_{12}-\vec{r}_1), \label{Stru-c-1}\\
  \vec{r}_{12} &=& a\vec{r}+b(\vec{r}_1-\vec{r})+c(\vec{r}_2-\vec{r}),  \label{Stru-c-2}\\
  \vec{r}_{22}-\vec{r}_2 &=& \gamma(\vec{r}_2-\vec{r})+\delta(\vec{r}_{12}-\vec{r}_2),\label{Stru-c-3}
\end{eqnarray}
where $\alpha,\beta,a,b,c,\gamma$ and $\delta$ are discrete functions from $\mathbb{Z}^2$ to $\mathbb{R}$.
\par The following proposition describes some centroaffine invariants included in the above structure equations. Exactly, all coefficients in Eqs. (\ref{Stru-c-1})-(\ref{Stru-c-3}) are centroaffine invariant.
\begin{prop}\label{Equ-pro}If two discrete centroaffine indefinite surface are centroaffinely equivalent, they have same discrete functions $a,b,c,\alpha,\beta,\gamma, \delta$.
\end{prop}
{\bf Proof.}
If two discrete centroaffine indefinite surface $\vec{r}(i,j)$ and $\vec{\bar{r}}(i,j)$ are centroaffinely equivalent, there exists a  non-degenerate matrix $P$ satisfying that
$\vec{r}(i,j)=P\vec{\bar{r}}(i,j).$  According to Eqs. (\ref{Stru-c-1})-(\ref{Stru-c-3}) and (\ref{ind-vec}), we obtain
\begin{gather}
  \alpha=\frac{\det[\vec{r}_1,\vec{r}_{11},\vec{r}_{12}]}{\det[\vec{r},\vec{r}_1,\vec{r}_{12}]}, \qquad
  \beta=\frac{\det[\vec{r},\vec{r}_{1},\vec{r}_{11}]}{\det[\vec{r},\vec{r}_1,\vec{r}_{12}]},  \\
  \gamma=\frac{\det[\vec{r}_2,\vec{r}_{22},\vec{r}_{12}]}{\det[\vec{r},\vec{r}_2,\vec{r}_{12}]}, \qquad
  \delta=\frac{\det[\vec{r},\vec{r}_{2},\vec{r}_{22}]}{\det[\vec{r},\vec{r}_2,\vec{r}_{12}]}, \\
  a=\frac{\det[\vec{r}_{12},\vec{r}_{1},\vec{r}_{2}]}{\det[\vec{r},\vec{r}_1,\vec{r}_{2}]},\quad
  b=\frac{\det[\vec{r},\vec{r}_{12},\vec{r}_{2}]}{\det[\vec{r},\vec{r}_1,\vec{r}_{2}]}, \quad
  c =\frac{\det[\vec{r},\vec{r}_{1},\vec{r}_{12}]}{\det[\vec{r},\vec{r}_1,\vec{r}_{2}]}.
\end{gather}
It can be easily seen that $$\bar{\alpha}=\alpha, \bar{\beta}=\beta, \bar{\gamma}=\gamma,\bar{\delta}=\delta, \bar{a}=a,\bar{b}=b,\bar{c}=c,$$
which show the discrete functions $a,b,c,\alpha,\beta,\gamma, \delta$ are invariant under centroaffine transformation.
\\ \rightline{$\Box$}
By Definition \ref{dics}, there are some limitation to the coefficients in Eqs. (\ref{Stru-c-1})-(\ref{Stru-c-3}).
\begin{rem}
\begin{itemize}
  \item[(1)] From Eq. (\ref{ind-vec}), that is, Condition (b) and Condition (c), we can derive that $a-b-c\neq0,\quad b\neq0,\quad c\neq0.$
  \item[(2)] $a=1$ means $\vec{r}_{12}-\vec{r}, \vec{r}_2-\vec{r}, \vec{r}_1-\vec{r}$ are coplanar, which implies the discrete centroaffine surface is a plane locally.
\end{itemize}
\end{rem}
Thus, in this paper we assume
\begin{equation}\label{assume}
  a-b-c\neq0,\quad bc\neq0,\quad a\neq1.
\end{equation}
Hence, the compatibility conditions of  Eqs. (\ref{Stru-c-1})-(\ref{Stru-c-3}) yield
\begin{eqnarray}
  \alpha &=& \frac{(1-a_1)(a-b-c)}{(a-1)b_1}, \label{Dis-Con-F} \\
  \gamma &=& \frac{(1-a_2)(a-b-c)}{(a-1)c_2}, \label{Dis-Con-S} \\
  \frac{c}{b} &=& \frac{c_{12}(a_2-b_2-c_2)}{b_{12}(a_1-b_1-c_1)},\\
  \frac{1-a_{12}}{b_{12}}&=& \frac{(a_2-1)(a_1-1)c}{(a-1)(a_2-b_2-c_2)}(1-K),
\end{eqnarray}
and
\begin{eqnarray}
  (a_2-1)b\beta_2+(a-1)(a_1-c_1) &=& (a-1)b_1\beta+(a_1-1)(a-c),  \\
  (a_1-1)c\delta_1+(a-1)(a_2-b_2) &=& (a-1)c_2\delta+(a_2-1)(a-b), \label{Dis-Con-L}
\end{eqnarray}
where $ K=\frac{[(a-1)b_1\beta+(a_1-c_1)(1-a)-(1-a_1)(a-c)][c_2\delta(a-1)+(a_2-1)(a-b)-(a-1)(a_2-b_2)]}{bc(a_2-1)(a_1-1)}.$
\begin{rem} In view of $a-b-c\neq0, a\neq1,$ Eqs. (\ref{Dis-Con-F}) and (\ref{Dis-Con-S}), it is obvious that
\begin{equation}
  \alpha\neq 0,\gamma\neq 0.
\end{equation}
\end{rem}
\par Integrable discrete versions of indefinite affine spheres have been constructed in \cite{Bobenko-3,Bo-Pin-99} by the following equations, which is shown that
the underlying discrete Gauss-Codazzi equations reduce to an integrable discrete Tzitzeica system.
\begin{thm}(The discrete Tzitzeica system). Discrete affine sphere are governed by the discrete Gauss equations
\begin{eqnarray}
  \vec{r}_{11}-\vec{r}_1 &=& \frac{H_1-1}{H_1(H-1)}(\vec{r}_1-r)+\frac{A}{H-1}(\vec{r}_{12}-\vec{r}_1),\label{AS-Con-F} \\
  \vec{r}_{12}+\vec{r} &=& H(\vec{r}_1+\vec{r}_2), \label{AS-Con-M}\\
  \vec{r}_{22}-\vec{r}_2 &=& \frac{H_2-1}{H_2(H-1)}(\vec{r}_2-\vec{r})+\frac{B}{H-1}(\vec{r}_{12}-\vec{r}_2). \label{AS-Con-L}
\end{eqnarray}
They are compatible modulo
\begin{gather}
  A_2=\frac{H_1}{H}A,\quad B_1=\frac{H_2}{H}B,\\
  H_{12}=\frac{H(H-1)}{H^2(H_1+H_2-H_1H_2)-H+ABH_1H_2},
\end{gather}
which is termed the discrete Tzitzeica system\cite{Bo-Pin-99}.
\end{thm}
By comparing Eqs. (\ref{Stru-c-1})-(\ref{Stru-c-3}) and Eqs. (\ref{AS-Con-F})-(\ref{AS-Con-L}), it is easy to see
\begin{prop}\label{BeqC} For a discrete centroaffine indefinite surface, if the coefficients satisfy that $b=c, a-b-c=-1$, it is a discrete affine sphere.
\end{prop}
Using Eqs. (\ref{Stru-c-1})-(\ref{Stru-c-3}), we get the following chain rule
\begin{align}
  \left[\vec{r}(m+1,n),\vec{r}(m+2,n),\vec{r}(m+1,n+1)\right]&\nonumber\\
   =&\left[\vec{r}(m,n),\vec{r}(m+1,n),\vec{r}(m,n+1)\right]A(m,n), \label{Tran-1}\\
 \left[\vec{r}(m,n+1),\vec{r}(m+1,n+1),\vec{r}(m,n+2)\right]& \nonumber\\
 =& \left[\vec{r}(m,n),\vec{r}(m+1,n),\vec{r}(m,n+1)\right]B(m,n)\label{Tran-2},
\end{align}
where
\begin{eqnarray}
  A(m,n) &=& \left(
          \begin{array}{ccc}
            0 & \beta_{mn}(a_{mn}-b_{mn}-c_{mn})-\alpha_{mn} & a_{mn}-b_{mn}-c_{mn} \\
            1 & 1+\alpha_{mn}+b_{mn}\beta_{mn}-\beta_{mn} & b_{mn} \\
            0 & c_{mn}\beta_{mn} & c_{mn} \\
          \end{array}
        \right),
   \\
  B(m,n) &=&  \left(
          \begin{array}{ccc}
            0 &  a_{mn}-b_{mn}-c_{mn} & \delta_{mn}(a_{mn}-b_{mn}-c_{mn})-\gamma_{mn}\\
            0 &  b_{mn}     &  b_{mn}\delta_{mn}\\
            1 &  c_{mn}     &  1+\gamma_{mn}+c_{mn}\delta_{mn}-\delta_{mn}\\
          \end{array}
        \right).
\end{eqnarray}
By a direct computation, we have
\begin{equation}\label{detA-B}
 |A|=c_{mn}\alpha_{mn},|B|=b_{mn}\gamma_{mn},
\end{equation}
where $|A|$ indicates the determinant of $A$.
\par Eqs. (\ref{Tran-1}) and (\ref{Tran-2}) give
\begin{align*}
        \left[\vec{r}(m+1,n+1),\vec{r}(m+2,n+1),\vec{r}(m+1,n+2)\right]& \\
        =\left[\vec{r}(m,n),\vec{r}(m+1,n),\vec{r}(m,n+1)\right]&A(m,n)B(m+1,n)  \\
         =\left[\vec{r}(m,n),\vec{r}(m+1,n),\vec{r}(m,n+1)\right]&B(m,n)A(m,n+1),
\end{align*}
which leads to
\begin{equation}\label{com-con}
 A(m,n)B(m+1,n) = B(m,n)A(m,n+1).
\end{equation}
\par Obviously, the result of this chain rule is accordance with Eqs. (\ref{Stru-c-1})-(\ref{Stru-c-3}), so they have same compatibility conditions. Thus, we conclude
\begin{prop}The equation $$ A(m,n)B(m+1,n) = B(m,n)A(m,n+1)$$ is equivalent to the compatibility conditions (\ref{Dis-Con-F})-(\ref{Dis-Con-L}).
\end{prop}

\begin{figure}[hbtp]
            \centering
           \includegraphics[width=.6\textwidth]{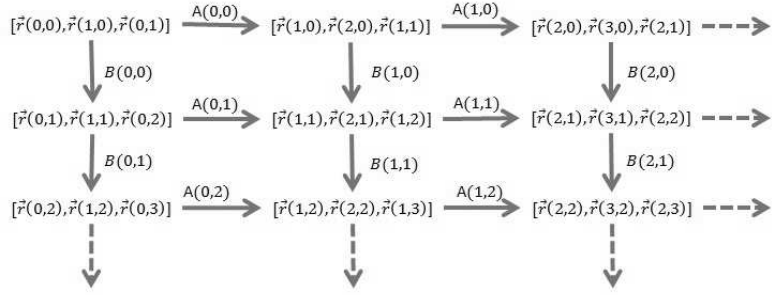}
            \caption{The chain graph of the points in proper order. }
            \label{fig-gpoint}
 \end{figure}

\begin{figure}[hbtp]
            \centering
           \includegraphics[width=.8\textwidth]{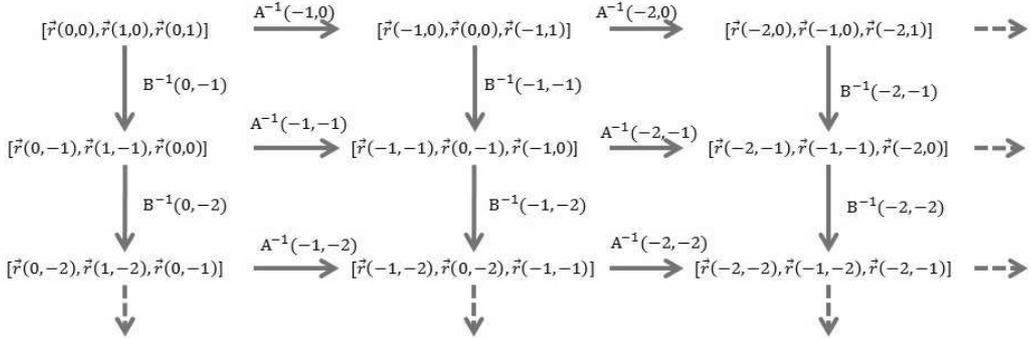}
            \caption{The chain graph of the points in reverse order. }
            \label{fig-gpoint1}
 \end{figure}

Given the initial three points $\vec{r}(0,0),\vec{r}(1,0),\vec{r}(0,1)$ with the discrete functions $a,b,c,\alpha,\beta,\gamma, \delta,$ all point groups in the discrete centroaffine surface can be generated according to the chain rule shown in Figure \ref{fig-gpoint} and Figure \ref{fig-gpoint1}. Especially, under a centroaffine transformation, the discrete functions $a,b,c,\alpha,\beta,\gamma, \delta,$ are invariant. On the other hand, there exists a centroaffine transformation which changes $\vec{r}(0,0),\vec{r}(1,0),\vec{r}(0,1)$ to $(1,0,0)^{\mathrm{Tran}}, (0,1,0)^{\mathrm{Tran}},(0,0,1)^{\mathrm{Tran}}.$ Hence we can always choose $$\vec{r}(0,0)=(1,0,0)^{\mathrm{Tran}}, \vec{r}(1,0)=(0,1,0)^{\mathrm{Tran}},\vec{r}(0,1)=(0,0,1)^{\mathrm{Tran}}.$$

\par Now, we can directly obtain the point group $\vec{r}(m,n),\vec{r}(m+1,n),\vec{r}(m,n+1)$ by the transition matrices, and it is similar for the point group $\vec{r}(-m,-n),\vec{r}(-m+1,-n),\vec{r}(-m,-n+1)$. By Figure \ref{fig-gpoint} and Figure \ref{fig-gpoint1},
\begin{align}
        \left[\vec{r}(m,n),\vec{r}(m+1,n),\vec{r}(m,n+1)\right]& \nonumber\\
        =\left[\vec{r}(0,0),\vec{r}(1,0),\vec{r}(0,1)\right]A(0,0)&A(1,0)\cdots A(m-1,0)B(m,0)B(m,1)\cdots B(m,n-1)\label{Transp}
\end{align}
and
\begin{align}
        \left[\vec{r}(-m,-n),\vec{r}(-m+1,-n),\vec{r}(-m,-n+1)\right]& \nonumber\\
        =\left[\vec{r}(0,0),\vec{r}(1,0),\vec{r}(0,1)\right]A^{-1}(-1,0)\cdots &A^{-1}(-m,0)B^{-1}(-m,-1)\cdots B^{-1}(-m,-n),\label{Transr}
\end{align}
where $m>0,n>0$ and $m,n\in\mathbb{Z}$.

\par Under a centroaffine transformation, the points $\vec{r}(0,0), \vec{r}(1,0),\vec{r}(0,1)$ can be changed to
 $$\vec{r}(0,0)=(1,0,0)^{\mathrm{Tran}}, \vec{r}(1,0)=(0,1,0)^{\mathrm{Tran}},\vec{r}(0,1)=(0,0,1)^{\mathrm{Tran}}.$$
 and the discrete functions $a,b,c,\alpha,\beta,\gamma, \delta$ are invariant. Thus the following proposition is obvious.
\begin{prop}\label{Exi-pro}Given the discrete function $a,b,c,\alpha,\beta,\gamma, \delta$ satisfying the relations (\ref{Dis-Con-F})-(\ref{Dis-Con-L}), there exists an unique discrete centroaffine indefinite surface under a centroaffine transformation.
\end{prop}
Then by Proposition \ref{Equ-pro} and Proposition \ref{Exi-pro}, we obtain
\begin{cor}\label{Lem-same}If two discrete centroaffine indefinite surface $\vec{r}(i,j)$ and $\vec{\bar{r}}(i,j)$ have same coefficients, that is, $$\bar{\alpha}(i,j)=\alpha(i,j), \bar{\beta}(i,j)=\beta(i,j), \bar{\gamma}(i,j)=\gamma(i,j),\bar{\delta}(i,j)=\delta(i,j),$$ $$\bar{a}(i,j)=a(i,j),\bar{b}(i,j)=b(i,j),\bar{c}(i,j)=c(i,j),$$ they are centroaffine equivalent.
\end{cor}
\par
\begin{figure}[hbtp]
            \centering
            \begin{tabular}{cc}
              \includegraphics[width=.45\textwidth]{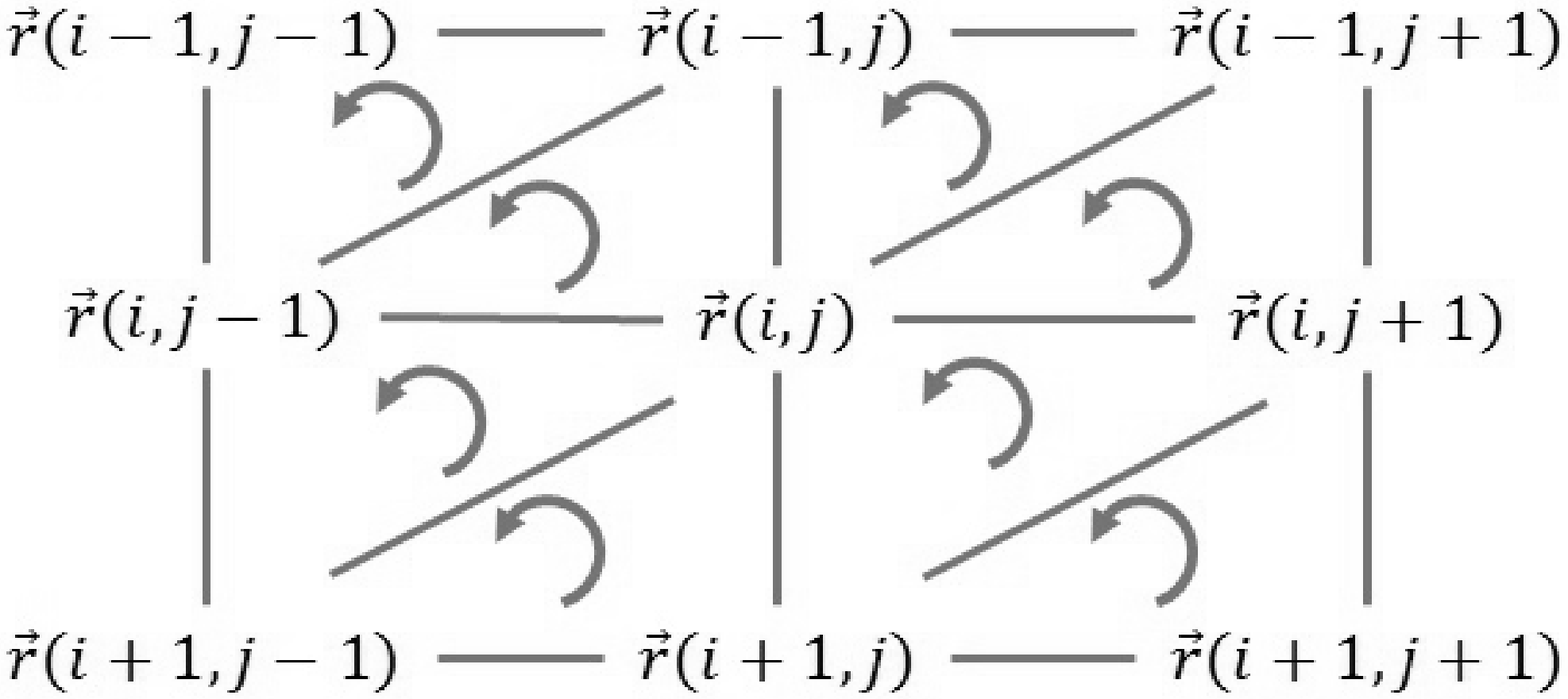} & \includegraphics[width=.45\textwidth]{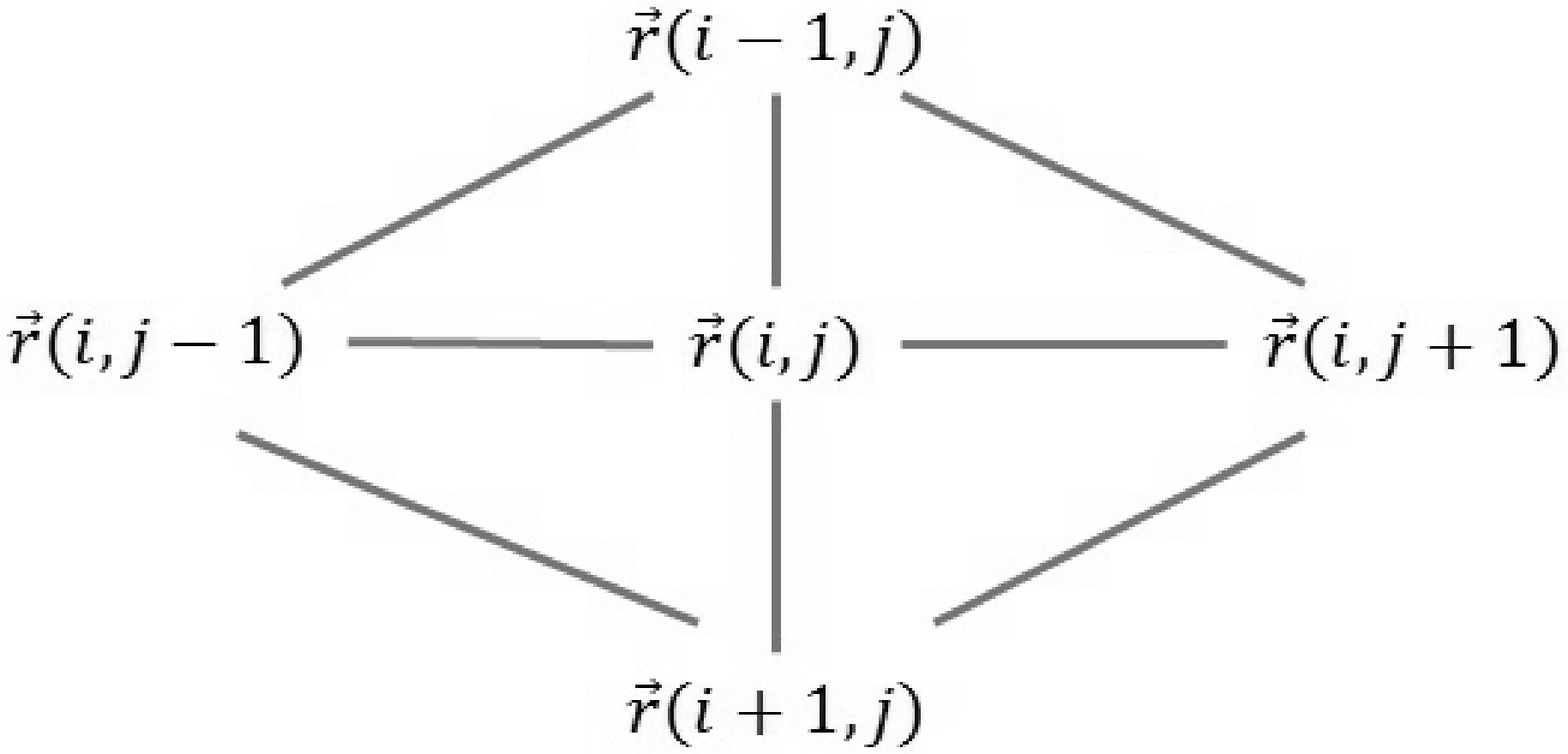}
            \end{tabular}
            \caption{{\it Left} : Each vertex neighborhood of a triangle mesh looks topologically like an oriented piece of the
plane. {\it Right} : The tangent plane of the point $\vec{r}(i,j)$. }
            \label{fig-DS}
 \end{figure}

 Given all point group $\{\vec{r}(i,j)\}$, we can generate the discrete surface by the following mode on the left side of Figure \ref{fig-DS}, and every triangle represents a face in the discrete centroaffine indefinite surface. Moreover, the plane $\pi$ can be regarded as {\it the tangent plane} at the point $\vec{r}(i,j)$ on the right side of Figure \ref{fig-DS}.
\par A discrete centroaffine indefinite surface $M=(V,E,F)$ is a simplicial 2-complex consisting of
  \begin{itemize}
    \item distinct vertices $$V=\{\vec{r}(i,j), i,j\in\mathbb{Z}\};$$
    \item oriented edges
    \begin{align*}
    E=&\{[\vec{r}(i,j-1),\vec{r}(i,j)],[\vec{r}(i,j),\vec{r}(i-1,j)],[\vec{r}(i-1,j+1),\vec{r}(i,j)],\\
    &[\vec{r}(i,j),\vec{r}(i,j+1)],[\vec{r}(i+1,j),\vec{r}(i,j)],[\vec{r}(i,j),\vec{r}(i+1,j-1)],i,j\in\mathbb{Z}\};
    \end{align*}
    \item oriented faces $$F=\{(\vec{r}(i,j),\vec{r}(i,j+1),\vec{r}(i+1,j)),(\vec{r}(i,j-1),\vec{r}(i+1,j-1),\vec{r}(i,j)),i,j\in\mathbb{Z}\}.$$
  \end{itemize}
Let star($\vec{r}(i,j)$) denote the triangles of $M$ that contain $\vec{r}(i,j)$ as a vertex. For an edge $\vec{e}$, let star($\vec{e}$) denote the (at most two) triangles of $M$ that contain $\vec{e}$ as an edge. Let $T=(\vec{p},\vec{q},\vec{r})$ denote an oriented triangle of $M$ with vertices $p,q,r\in M$. The volume of an oriented surface $M$ is the oriented volume enclosed by the cone of the surface over the origin in $\R^3$:
\begin{equation}
  \textbf{Vol}\{M\}:=\frac{1}{6}\sum_{T=(\vec{p},~\vec{q},~\vec{r})\in M}\det[\vec{p},\vec{q},\vec{r}].
\end{equation}
By simplification, we take the following notations:
$$d(i,j)=a(i,j)-b(i,j)-c(i,j),\quad \textbf{V}(i,j)=\textbf{Vol}\{(\vec{r}(i,j),\vec{r}(i+1,j),\vec{r}(i,j+1))\}.$$
Thus
\begin{align}
 \textbf{Vol}\{\mathrm{star}([\vec{r}(i,j),\vec{r}(i-1,j+1)])\}&=[1-d(i-1,j)]\textbf{V}(i-1,j).\\
  \textbf{Vol}\{(\vec{r}(i-1,j),\vec{r}(i,j),\vec{r}(i,j+1))\}&=c(i-1,j)\textbf{V}(i-1,j).
\end{align}
 We also have
 \begin{align}
 \textbf{Vol}\{\mathrm{star}([\vec{r}(i,j),\vec{r}(i+1,j-1)])\}&=[1-d(i,j-1)]\textbf{V}(i,j-1).\\
  \textbf{Vol}\{(\vec{r}(i,j-1),\vec{r}(i+1,j),\vec{r}(i,j))\}&=b(i,j-1)\textbf{V}(i,j-1).
\end{align}
On the other hand
 \begin{align}
 \textbf{Vol}\{(\vec{r}(i,j-1),\vec{r}(i,j),\vec{r}(i-1,j))\}=-d(i-1,j-1)\textbf{V}(i-1,j-1).
\end{align}
Finally we get
 \begin{align}
 \textbf{Vol}\{\mathrm{star}(\vec{r}(i,j))\}&=-d(i-1,j-1)\textbf{V}(i-1,j-1)+[1-d(i-1,j)]\textbf{V}(i-1,j)\nonumber\\
 &+[1-d(i,j-1)]\textbf{V}(i,j-1)+\textbf{V}(i,j)\nonumber\\
 &=[-d_{i-1,j-1}+(1-d_{i-1,j})|B_{i-1,j-1}|+(1-d_{i,j-1})|A_{i-1,j-1}|\nonumber\\
 &+|A_{i-1,j-1}B_{i,j-1}|~]\textbf{V}_{i-1,j-1}\nonumber\\
 &=[-d_{i-1,j-1}+(1-d_{i-1,j})|B_{i-1,j-1}|+(1-d_{i,j-1})|A_{i-1,j-1}|\nonumber\\
 &+|A_{i-1,j-1}B_{i,j-1}|~]\times|A(0,0)A(1,0)\cdots A(i-2,0)|\nonumber\\
 &\times|B(i-1,0)\cdots B(i-1,j-2)|V(0,0).\label{Vol-rij}
\end{align}
In the tangent plane,
\begin{align}
 \textbf{Vol}_1\{\mathrm{star}(\vec{r}(i,j))\}&=-d(i-1,j-1)\textbf{V}(i-1,j-1)+c(i-1,j)\textbf{V}(i-1,j)\nonumber\\
 &+b(i,j-1)\textbf{V}(i,j-1)+\textbf{V}(i,j)\nonumber\\
 &=[-d_{i-1,j-1}+c_{i-1,j}|B_{i-1,j-1}|+b_{i,j-1}|A_{i-1,j-1}|\nonumber\\
 &+|A_{i-1,j-1}B_{i,j-1}|~]\textbf{V}_{i-1,j-1}\nonumber\\
 &=[-d_{i-1,j-1}+c_{i-1,j}|B_{i-1,j-1}|+b_{i,j-1}|A_{i-1,j-1}|\nonumber\\
 &+|A_{i-1,j-1}B_{i,j-1}|~]\times|A(0,0)A(1,0)\cdots A(i-2,0)|\nonumber\\
 &\times|B(i-1,0)\cdots B(i-1,j-2)|V(0,0).\label{tVol-rij}
\end{align}
\section{Laplacian operator on centroaffine simplicial surface.}
Let $f:\mathcal{S}\rightarrow\R^3$ be a simplicial surface $S:=f(\mathcal{S})$. Then the half square edge energy of the simplicial surface $S$  is given by
$$E(S)=\frac{1}{2}\sum_{(i,j)\in E}||f(i)-f(j)||^2,$$
where $$||f(i)-f(j)||^2=\sum_{k=1}^3|f_k(i)-f_k(j)|^2$$
and $f_k$ is the $k-$th coordinate function. Then
$$E(S)=\sum_{k=1}^3E(f_k),$$
where $E(f_k)$ is the Dirichlet energy of the $k-$th coordinate function.
Using the similar method as in \cite{Bo-Pin-07}, the gradient of $E(S)$ at the vertex $f(i)$ is equal to
$$\nabla_{f(i)}E(S)=\sum_{j:(i,j)\in E}(f(i)-f(j)).$$
In analogy to the smooth case, the Laplacian operator is defined as the
gradient of the Dirichlet energy. So we obtain the {\it combinatorial Laplacian operator}
\begin{equation}\label{Laplacian}
 (\Delta f)(i)=\sum_{j:(i,j)\in E}(f(i)-f(j)).
\end{equation}
This Laplacian operator is very useful due to its simple definition using only the connectivity of the mesh.  Besides that matrix is symmetric, and so it has real eigenvalues and orthogonal eigenvectors. On the other hand, from the following proposition, it is closely related to centroaffine transformation.
\begin{prop}
$\Delta f=0$ and $\Delta f= sf$ is centroaffine invariant,  where $s$ is a centroaffine invariant.
\end{prop}
{\bf Proof.} Under a centroaffine transformation $\tilde{f}=Pf$, where $P$ is a non-degenerate matrix, from Eq. (\ref{Laplacian}) we can show that
\begin{eqnarray*}
  (\Delta \tilde{f})(i)&=&\sum_{j:(i,j)\in E}(\tilde{f}(i)-\tilde{f}(j)) \\
   &=& \sum_{j:(i,j)\in E}(Pf(i)-Pf(j)) \\
   &=&  P\sum_{j:(i,j)\in E}(f(i)-f(j)).
\end{eqnarray*}
Hence, $\Delta f=0$ and $\Delta f= sf$ can generate $\Delta \tilde{f}=0$ and $\Delta \tilde{f}= s\tilde{f}$ respectively.
\\ \rightline{$\Box$}
For a discrete centroaffine indefinite surface $S=(V,E,F)$, by Eq.(\ref{Laplacian}) and Figure \ref{fig-DS}, we have
\begin{equation}\label{Lap1}
  \Delta \vec{r}(i,j)=6\vec{r}(i,j)-\vec{r}(i-1,j)-
\vec{r}(i-1,j+1)-\vec{r}(i,j-1)-\vec{r}(i,j+1)-\vec{r}(i+1,j-1)-\vec{r}(i+1,j).
\end{equation}
\begin{defn}
A discrete centroaffine indefinite surface is called harmonic if $\Delta \vec{r}=0.$
\end{defn}

A discrete centroaffine indefinite surface $S=(V,E,F)$ is harmonic if and only if the point $\vec{r}(i,j)$ lies the center of gravity of its immediate neighbors $\vec{r}(i-1,j)$,
$\vec{r}(i-1,j+1)$,$\vec{r}(i,j-1)$,$\vec{r}(i,j+1)$,$\vec{r}(i+1,j-1)$ and $\vec{r}(i+1,j)$, distributing the vertices over the space in a good way. Now, the following two corollaries are obvious.
\begin{cor}
A discrete centroaffine indefinite surface is harmonic if and only if
$$\vec{r}(i,j)=\frac{\vec{r}(i-1,j)+
\vec{r}(i-1,j+1)+\vec{r}(i,j-1)+\vec{r}(i,j+1)+\vec{r}(i+1,j-1)+\vec{r}(i+1,j)}{6}.$$
\end{cor}
\begin{cor}A discrete centroaffine indefinite surface satisfies that $\Delta \vec{r}= s\vec{r}$ if and only if
$(6+s)\vec{r}(i,j)= \vec{r}(i-1,j)+
\vec{r}(i-1,j+1)+\vec{r}(i,j-1)+\vec{r}(i,j+1)+\vec{r}(i+1,j-1)+\vec{r}(i+1,j).$
\end{cor}
\section{Discrete centroaffine indefinite surface with constant coefficients.}
In this section, we will consider the discrete centroaffine indefinite surface with constant coefficients in structure equations (\ref{Stru-c-1})-(\ref{Stru-c-3}).
If all the coefficients are constant, the compatibility conditions (\ref{Dis-Con-F})-(\ref{Dis-Con-L}) may be written as
\begin{gather}\label{com-con-c}
  -\alpha b=-\gamma c=a-b-c=bc(\beta\delta-1).
\end{gather}
The transition equations (\ref{Transp}) and (\ref{Transr}) are changed to
\begin{gather}
        \left[\vec{r}(m,n),\vec{r}(m+1,n),\vec{r}(m,n+1)\right]=\left[\vec{r}(0,0),\vec{r}(1,0),\vec{r}(0,1)\right]A^mB^n,\label{Ite}\\
        \left[\vec{r}(-m,-n),\vec{r}(-m+1,-n),\vec{r}(-m,-n+1)\right]=\left[\vec{r}(0,0),\vec{r}(1,0),\vec{r}(0,1)\right](A^{-1})^m(B^{-1})^n,
\end{gather}
where
\begin{eqnarray}
  A &=& \left(
          \begin{array}{ccc}
            0 & \beta(a-b-c)-\alpha & a-b-c \\
            1 & 1+\alpha+b\beta-\beta & b \\
            0 & c\beta & c \\
          \end{array}
        \right),
   \\
  B &=&  \left(
          \begin{array}{ccc}
            0 &  a-b-c & \delta(a-b-c)-\gamma\\
            0 &  b     &  b\delta\\
            1 &  c     &  1+\gamma+c\delta-\delta\\
          \end{array}
        \right),
\end{eqnarray}
 $m>0,n>0$ and $m,n\in\mathbb{Z}$. As in Eqs. (\ref{detA-B}) and (\ref{com-con}), we have
 \begin{equation}
 |A|=c\alpha,\quad|B|=b\gamma,\quad  AB=BA.
\end{equation}
In Proposition \ref{BeqC}, we obtain the relation between a discrete centroaffine indefinite surface and a discrete affine sphere. The following proposition gives another result.
\begin{prop}\label{Pro-Aff-sph}A discrete centroaffine indefinite surface with constant coefficients is a discrete affine sphere if and only if $|A|=|B|=1$.
\end{prop}
{\bf Proof.} Firstly, if $|A|=|B|=1$, which implies $\alpha=\frac{1}{c}, \gamma=\frac{1}{b}$, Eq. (\ref{com-con-c}) generates
$\frac{b}{c}=\frac{c}{b}.$ Then we have $b=-c$ or $b=c$. If $b=-c$, it is easy to see $\alpha=-\frac{1}{b}$ and $a=a-b-c=-\alpha b=1$, which is contradictory to the assumption $a\neq1$ in Eq. (\ref{assume}).  If $b=c$, we get $\alpha=\frac{1}{b}$ and $a-b-c=-\alpha b=-1$, which satisfies the condition of a discrete affine sphere in Proposition \ref{BeqC}.
\par On the other hand, if $b=c, a-b-c=-1$, Eq. (\ref{com-con-c}) gives $\alpha=\gamma=\frac{1}{b}=\frac{1}{c}$. It immediately shows that $|A|=c\alpha=1,\quad|B|=b\gamma=1$.
\\ \rightline{$\Box$}
The following proposition shows an interest property of a discrete centroaffine indefinite surface with the constant coefficients.
\begin{prop}\label{Prop-locally}
A discrete centroaffine indefinite surface with the constant coefficients is self-equivalent locally, that is, the patch $$\vec{r}(i,j),\cdots,\vec{r}(i+m,j+n)$$ is centroaffine equivalent to the patch $$\vec{r}(k,l),\cdots,\vec{r}(k+m,l+n),$$ where $i,j,k,l,m,n\in\mathbb{Z}$
\end{prop}
{\bf Proof.}
For the vector groups $\{\vec{r}(i,j),\vec{r}(i+1,j),\vec{r}(i,j+1)\}$ and $\{\vec{\bar{r}}(k,l),\vec{\bar{r}}(k+1,l),\vec{\bar{r}}(k,l+1)\}$,  there exists a non-degenerate matrix $P$ satisfying that $$\{\vec{r}(i,j),\vec{r}(i+1,j),\vec{r}(i,j+1)\}=P\{\vec{\bar{r}}(k,l),\vec{\bar{r}}(k+1,l),\vec{\bar{r}}(k,l+1)\}.$$
 Since the functions $\alpha,\beta,a,b,c,\gamma,\delta$ are constant, by Eqs. (\ref{Stru-c-1})-(\ref{Stru-c-3}), we get
 $$\{\vec{r}(i,j),\cdots,\vec{r}(i+m,j+n)\}=P\{\vec{r}(k,l),\cdots,\vec{r}(k+m,l+n)\},$$
  which completes the proof.
\\ \rightline{$\Box$}
\begin{rem}From the proposition \ref{Prop-locally}, locally, we can arbitrarily choose a patch to display a discrete centroaffine indefinite surface with constant coefficients.
\end{rem}
Since the coefficients $a,b,c,\alpha,\beta,\gamma,\delta$ and $\Delta\vec{r}=0$ are invariant under centroaffine transformation, we can choose $$\vec{r}(0,0)=(1,0,0)^{\mathrm{Tran}}, \vec{r}(1,0)=(0,1,0)^{\mathrm{Tran}},\vec{r}(0,1)=(0,0,1)^{\mathrm{Tran}}.$$
to obtain the following proposition by a direct computation.
\begin{prop}
A discrete centroaffine indefinite surface $S=(V,E,F)$ with constant coefficients is harmonic if and only if it satisfies that
\begin{gather}
 (1+\alpha-\beta)(\frac{1}{\alpha}+\frac{b}{c})+(1+\gamma-\delta)(\frac{1}{\gamma}+\frac{c}{b})-\frac{c}{b}-\frac{b}{c}=6,\label{Harmonic-F} \\
 \frac{\delta}{\gamma}(1+\alpha)=\frac{1}{\alpha}(1+\gamma)-\frac{1}{b}-1,\\
 \frac{\beta}{\alpha}(1+\gamma)=\frac{1}{\gamma}(1+\alpha)-\frac{1}{c}-1.\label{Harmonic-L}
\end{gather}
\end{prop}
\par Finally, according to above results, we give some examples for the discrete centroaffine indefinite surface with constant coefficients.
\par
\begin{figure}[hbtp]
            \centering
            \begin{tabular}{c}
              \includegraphics[width=.4\textwidth]{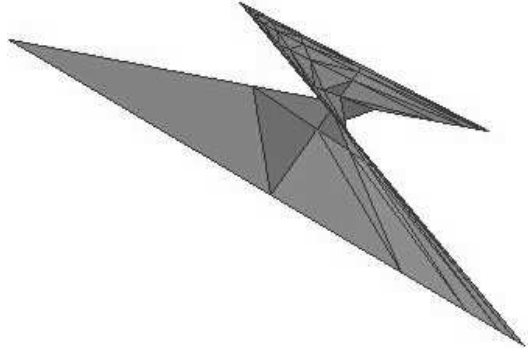}
            \end{tabular}
            \caption{A discrete centroaffine indifinite surface analogy of saddle surface. }
            \label{figure-saddle}
 \end{figure}

 {\bf Example 1.} It is not hard to check that $\beta=\delta=0, \alpha=\gamma=b=c=0.5$ and $a=0.75$ satisfy Eq. (\ref{com-con-c}), and we get a discrete surface like saddle surface shown in Figure \ref{figure-saddle}.
\par
\begin{figure}[hbtp]
            \centering
            \begin{tabular}{cc}
              \includegraphics[width=.25\textwidth]{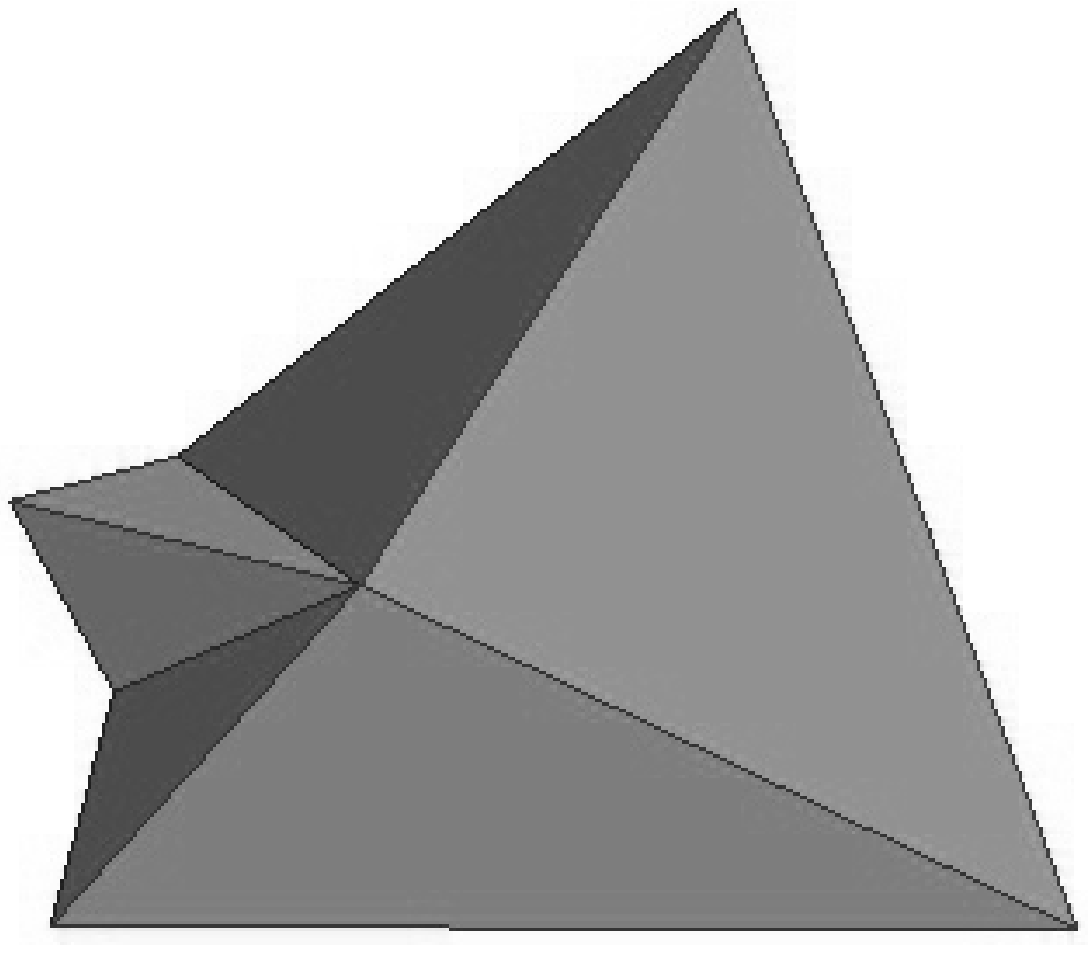} & \includegraphics[width=.3\textwidth]{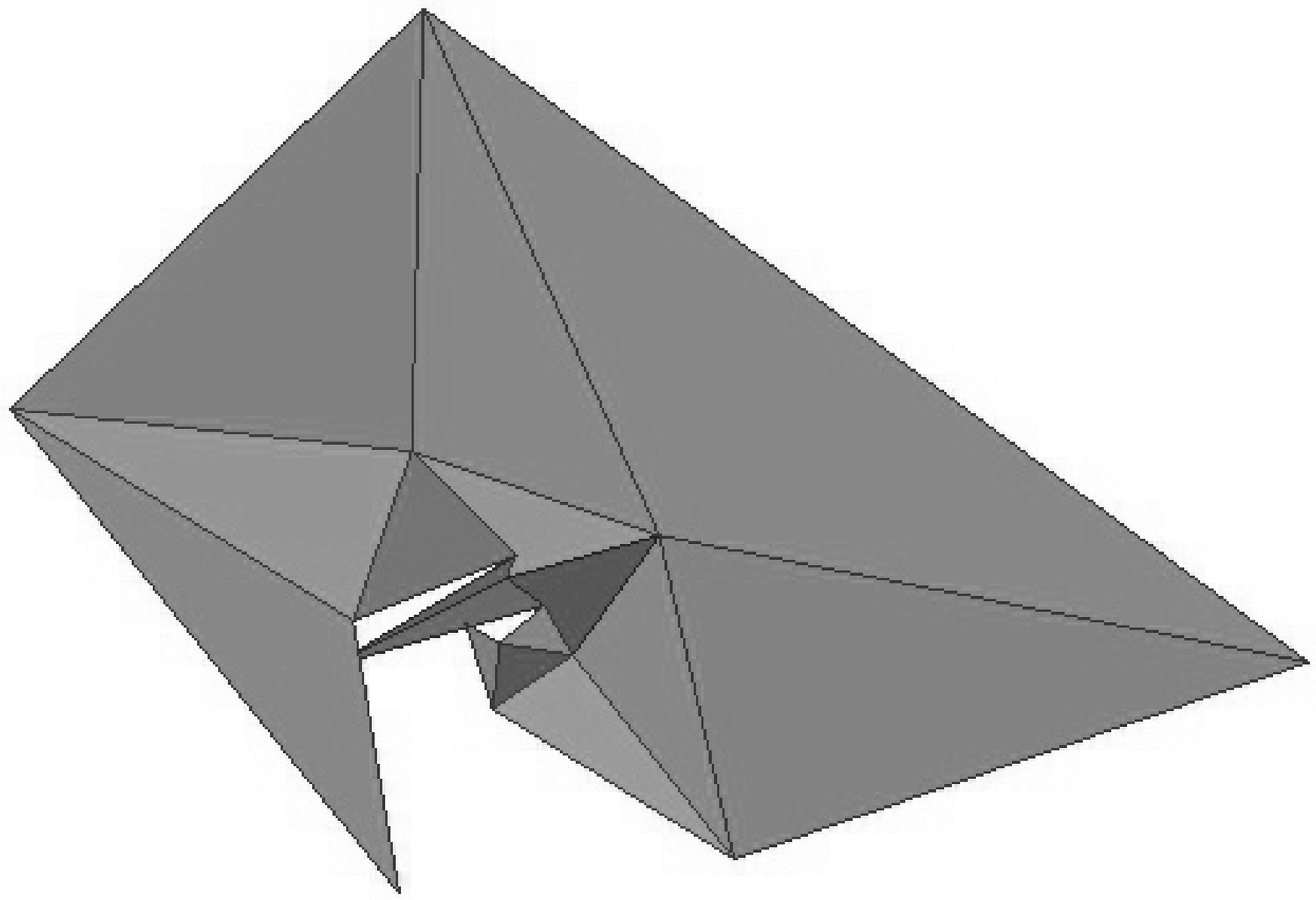}
            \end{tabular}
            \caption{Graph of a harmonic discrete centroaffine indefinite surface. }
            \label{figure-exp1}
 \end{figure}

 {\bf Example 2.} From the compatibility conditions (\ref{com-con-c}) and Eqs. (\ref{Harmonic-F})-(\ref{Harmonic-L}), by assuming $\alpha=1$, we get a harmonic discrete centroaffine indefinite surface with
$$a=-b=c=-\frac{1}{3}, \beta=\delta=2, \gamma=-1.$$
In particular, we obtain the coordinates of the corresponding points
\begin{gather*}
\vec{r}(0,0)=(1,0,0)^{\mathrm{Tran}},\vec{r}(1,0)=(0,1,0)^{\mathrm{Tran}},\vec{r}(0,1)=(0,0,1)^{\mathrm{Tran}},\\
\vec{r}(-1,1)=(1,1,-5)^{\mathrm{Tran}},\vec{r}(0,-1)=(2,-2,1)^{\mathrm{Tran}},\vec{r}(1,-1)=(3,1,1)^{\mathrm{Tran}},\\
\vec{r}(-1,0)=(0,-1,2)^{\mathrm{Tran}},\vec{r}(-1,-1)=(-1,-1,-1)^{\mathrm{Tran}}.
\end{gather*}
By these points, we obtain a local graph of a harmonic discrete centroaffine indefinite surface in the left of Figure \ref{figure-exp1}. With the iterative formula (\ref{Ite}), we generate 16 points and display it as a graph of a harmonic discrete centroaffine indefinite surface in the right of Figure \ref{figure-exp1}.
\par
\begin{figure}[hbtp]
            \centering
            \begin{tabular}{ccc}
               \includegraphics[width=.2\textwidth]{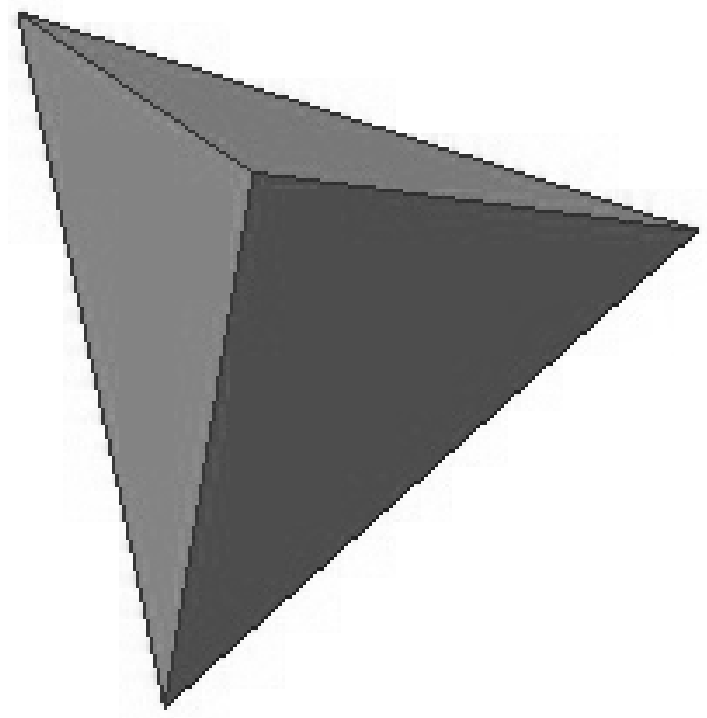}&\includegraphics[width=.23\textwidth]{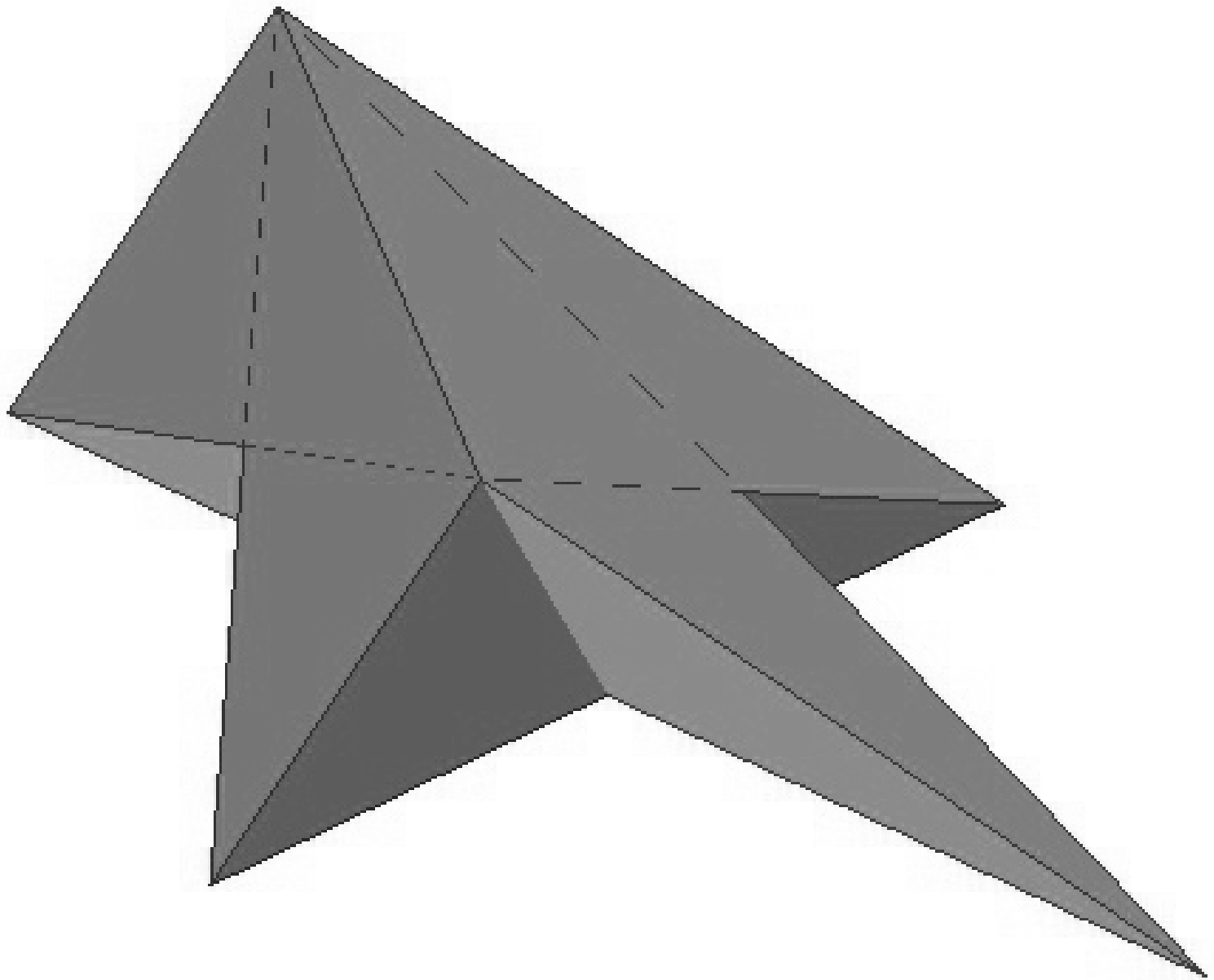}& \includegraphics[width=.13\textwidth]{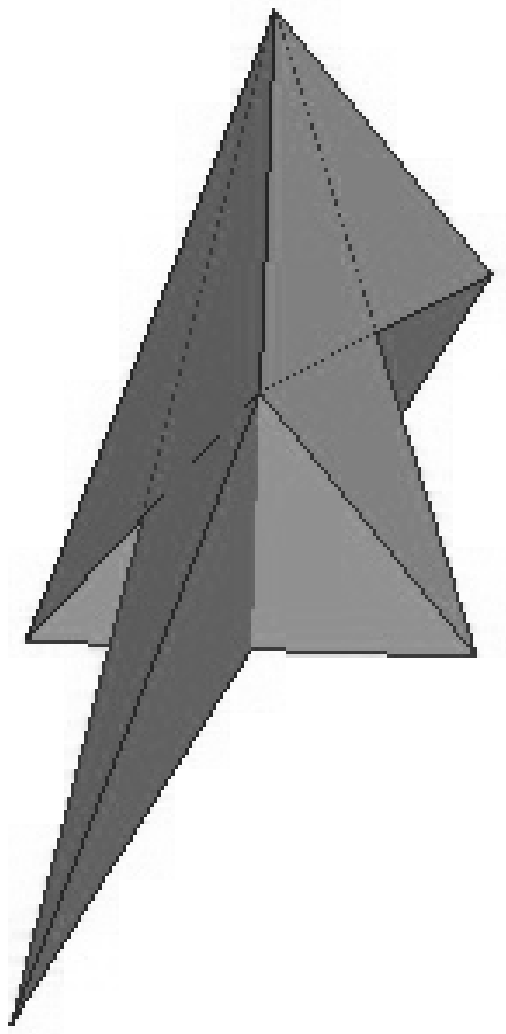}
            \end{tabular}
            \caption{Discrete affine sphere with $\Delta\vec{r}=8\vec{r}$. }
            \label{figure-AS}
 \end{figure}

 {\bf Example 3.} In the smooth case, a centroaffine surface with $\Delta_gx=-2x$ is a centroaffine minimal affine sphere. Here we consider a discrete affine sphere of constant coefficients with $\Delta \vec{r}=s\vec{r}$, where $s$ is constant. Using Eqs. (\ref{Lap1}), (\ref{com-con-c}) and Proposition \ref{Pro-Aff-sph}, by a simple calculation, we get $s=8$. Furthermore, we can obtain $a=-3,b=c=\alpha=\gamma=-1,\beta\delta=0$.
Especially, if $\beta=\delta=0$, the discrete affine sphere is the face of a tetrahedron with $\vec{r}_{11}=\vec{r}_{22}=\vec{r}$ and $\vec{r}_{12}=-\vec{r}-\vec{r}_1-\vec{r}_2$, as shown in the left of Figure \ref{figure-AS}.
If $\beta=0,\delta\neq0$,  the structure equations of the discrete affine sphere are
$$\vec{r}_{11}=\vec{r},\quad \vec{r}_{12}=-\vec{r}-\vec{r}_1-\vec{r}_2,\quad \vec{r}_{22}=(1-\delta)\vec{r}-\delta\vec{r}_1-2\delta\vec{r}_2.$$ Using $\delta=1$ and $\delta=-1$, we get the graphs shown with 9 vertices in the middle of Figure \ref{figure-AS} and in the right of Figure \ref{figure-AS}, respectively. In fact, there are only 6 vertices according to $\vec{r}_{11}=\vec{r}$.
\par From Proposition \ref{BeqC}, Proposition \ref{Pro-Aff-sph}, Eqs. (\ref{Vol-rij}) and (\ref{tVol-rij}) we have
\begin{cor} For a discrete affine sphere with constant coefficients we have
$$\textbf{V}(i,j)=\textbf{V}(0,0),\quad\frac{\textbf{Vol}\{\mathrm{star}(\vec{r}(i,j))\}}{\textbf{V}(0,0)}=6,\quad\frac{\textbf{Vol}_1\{\mathrm{star}(\vec{r}(i,j))\}}{\textbf{V}(0,0)}=2+b+c.$$
\end{cor}

\section{Locally convexity of a discrete centroaffine indefinite surface.}
Obviously, a discrete centroaffine indefinite surface is convex at the point $\vec{r}(n_1,n_2)$  if these points
$$\vec{r}_{\bar{1}\bar{1}}, \vec{r}_{\bar{1}\bar{2}}, \vec{r}_{\bar{2}\bar{2}}, \vec{r}_{\bar{1}2}, \vec{r}_{1\bar{2}}, \vec{r}_{11}, \vec{r}_{12}, \vec{r}_{22}$$ lie on the same side of the tangent plane $\pi$ of $\vec{r}(n_1,n_2)$, which implies the following determinants have the same sign.
\begin{gather}
\begin{split}
&[\vec{r}_1-\vec{r}, \vec{r}_2-\vec{r},\vec{r}_{\bar{1}\bar{1}}-\vec{r}],\quad [\vec{r}_1-\vec{r}, \vec{r}_2-\vec{r},\vec{r}_{\bar{1}\bar{2}}-\vec{r}],  \quad [\vec{r}_1-\vec{r}, \vec{r}_2-\vec{r},\vec{r}_{\bar{2}\bar{2}}-\vec{r}], \\
&[\vec{r}_1-\vec{r}, \vec{r}_2-\vec{r},\vec{r}_{1\bar{2}}-\vec{r}], \quad[\vec{r}_1-\vec{r}, \vec{r}_2-\vec{r},\vec{r}_{\bar{1}2}-\vec{r}],\quad
[\vec{r}_1-\vec{r}, \vec{r}_2-\vec{r},\vec{r}_{11}-\vec{r}],  \\
&[\vec{r}_1-\vec{r}, \vec{r}_2-\vec{r},\vec{r}_{12}-\vec{r}], \quad [\vec{r}_1-\vec{r}, \vec{r}_2-\vec{r},\vec{r}_{22}-\vec{r}].
\end{split}\label{Det-all}
\end{gather}
Firstly, from Eq. (\ref{Stru-c-2}), it is easy to check
\begin{equation}
\begin{split}
[\vec{r}_1-\vec{r}, \vec{r}_2-\vec{r},\vec{r}_{12}-\vec{r}]&=[\vec{r}_1-\vec{r}, \vec{r}_2-\vec{r},(a-1)\vec{r}]\\
&=(a-1)[\vec{r},\vec{r}_1,\vec{r}_2].
\end{split}\label{Det1}
\end{equation}
On the other hand, by Eq. (\ref{Stru-c-3}) we have
\begin{equation}
\begin{split}
[\vec{r}_1-\vec{r}, \vec{r}_2-\vec{r},\vec{r}_{22}-\vec{r}]&=[\vec{r}_1-\vec{r}, \vec{r}_2-\vec{r},\vec{r}_{22}-\vec{r}_2+\vec{r}_2-\vec{r}]\\
&=[\vec{r}_1-\vec{r}, \vec{r}_2-\vec{r},\delta(\vec{r}_{12}-\vec{r}_2)]\\
&=(a-1)\delta[\vec{r},\vec{r}_1,\vec{r}_2].
\end{split}\label{Det2}
\end{equation}
Since $a-1\neq 0$ by the assumption (\ref{assume}) ,  comparison of Eq. (\ref{Det1}) and Eq. (\ref{Det2}) shows that
\begin{equation}
  \delta\geq 0.
\end{equation}
By Eq. (\ref{Stru-c-1}), we may examine
\begin{equation}
\begin{split}
[\vec{r}_1-\vec{r}, \vec{r}_2-\vec{r},\vec{r}_{11}-\vec{r}]&=[\vec{r}_1-\vec{r}, \vec{r}_2-\vec{r},\vec{r}_{11}-\vec{r}_1+\vec{r}_1-\vec{r}]\\
&=[\vec{r}_1-\vec{r}, \vec{r}_2-\vec{r},\beta(\vec{r}_{12}-\vec{r}_2)]\\
&=(a-1)\beta[\vec{r},\vec{r}_1,\vec{r}_2].
\end{split}
\end{equation}
Similarly, it is convenient to conclude
\begin{equation}
  \beta\geq0.
\end{equation}
Again, from Eqs. (\ref{Stru-c-1})-(\ref{Stru-c-3}), it is clear that
\begin{equation}
\begin{split}
[\vec{r}_1-\vec{r}, \vec{r}_2-\vec{r},\vec{r}_{\bar{1}2}-\vec{r}]&=[\vec{r}_1-\vec{r}, \vec{r}_2-\vec{r},\vec{r}_{\bar{1}2}-\vec{r}_{\bar{1}}+\vec{r}_{\bar{1}}-\vec{r}]\\
&=[\vec{r}_1-\vec{r}, \vec{r}_2-\vec{r},\frac{1}{c_{\bar{1}}}(\vec{r}_2-a_{\bar{1}}\vec{r}_{\bar{1}})]\\
&=[\vec{r}_1-\vec{r}, \vec{r}_2-\vec{r},\frac{1}{c_{\bar{1}}}(\vec{r}_2-a_{\bar{1}}(\vec{r}_{\bar{1}}-\vec{r}+\vec{r}))]\\
&=-\frac{1}{c_{\bar{1}}}(a_{\bar{1}}-1)[\vec{r},\vec{r}_1,\vec{r}_2].
\end{split}
\end{equation}
By the assumption (\ref{assume}), comparing the above equation with Eq. (\ref{Det1}) leads to
\begin{equation}
  -\frac{1}{c_{\bar{1}}}(a_{\bar{1}}-1)(a-1)>0.
\end{equation}
In like manner we shall derive
\begin{equation}
\begin{split}
[\vec{r}_1-\vec{r}, \vec{r}_2-\vec{r},\vec{r}_{\bar{2}1}-\vec{r}]&=[\vec{r}_1-\vec{r}, \vec{r}_2-\vec{r},\vec{r}_{\bar{2}1}-\vec{r}_{\bar{2}}+\vec{r}_{\bar{2}}-\vec{r}]\\
&=[\vec{r}_1-\vec{r}, \vec{r}_2-\vec{r},\frac{1}{b_{\bar{2}}}(\vec{r}_1-a_{\bar{2}}\vec{r}_{\bar{2}})]\\
&=[\vec{r}_1-\vec{r}, \vec{r}_2-\vec{r},\frac{1}{b_{\bar{1}}}(\vec{r}_1-a_{\bar{2}}(\vec{r}_{\bar{2}}-\vec{r}+\vec{r}))]\\
&=-\frac{1}{b_{\bar{1}}}(a_{\bar{2}}-1)[\vec{r},\vec{r}_1,\vec{r}_2],
\end{split}
\end{equation}
and
\begin{equation}
  -\frac{1}{b_{\bar{2}}}(a_{\bar{2}}-1)(a-1)>0.
\end{equation}
Also, we have
\begin{equation}
\begin{split}
[\vec{r}_1-\vec{r}, \vec{r}_2-\vec{r},\vec{r}_{\bar{2}\bar{2}}-\vec{r}]&=[\vec{r}_1-\vec{r}, \vec{r}_2-\vec{r},\vec{r}_{\bar{2}\bar{2}}-\vec{r}_{\bar{2}}+\vec{r}_{\bar{2}}-\vec{r}]\\
&=[\vec{r}_1-\vec{r}, \vec{r}_2-\vec{r},\frac{\delta_{\bar{2}\bar{2}}}{\gamma_{\bar{2}\bar{2}}}(\vec{r}_{\bar{2}1}-\vec{r}_{\bar{2}})]\\
&=[\vec{r}_1-\vec{r}, \vec{r}_2-\vec{r},\frac{\delta_{\bar{2}\bar{2}}}{\gamma_{\bar{2}\bar{2}}}(\vec{r}_{\bar{2}1}-\vec{r}+\vec{r}-\vec{r}_{\bar{2}})]\\
&=-\frac{1}{b_{\bar{1}}}(a_{\bar{2}}-1)\frac{\delta_{\bar{2}\bar{2}}}{\gamma_{\bar{2}\bar{2}}}[\vec{r},\vec{r}_1,\vec{r}_2],
\end{split}
\end{equation}
and
\begin{equation}
  \frac{\delta_{\bar{2}\bar{2}}}{\gamma_{\bar{2}\bar{2}}}\geq0.
\end{equation}
Similarly,
\begin{equation}
\begin{split}
[\vec{r}_1-\vec{r}, \vec{r}_2-\vec{r},\vec{r}_{\bar{1}\bar{2}}-\vec{r}]
&=[\vec{r}_1-\vec{r}, \vec{r}_2-\vec{r},\vec{r}_{\bar{1}\bar{2}}-\vec{r}_{\bar{1}}+\vec{r}_{\bar{1}}-\vec{r}]\\
&=[\vec{r}_1-\vec{r}, \vec{r}_2-\vec{r},-\frac{1}{\gamma_{\bar{1}\bar{2}}}(\vec{r}_{\bar{1}2}-\vec{r}_{\bar{1}})]\\
&=[\vec{r}_1-\vec{r}, \vec{r}_2-\vec{r},-\frac{1}{\gamma_{\bar{1}\bar{2}}}(\vec{r}_{\bar{1}2}-\vec{r}+\vec{r}-\vec{r}_{\bar{1}})]\\
&=-\frac{1}{c_{\bar{1}}}(a_{\bar{1}}-1)\frac{-1}{\gamma_{\bar{1}\bar{2}}}[\vec{r},\vec{r}_1,\vec{r}_2],
\end{split}
\end{equation}
and
\begin{equation}
 \gamma_{\bar{1}\bar{2}}<0.
\end{equation}
Noticing that
\begin{equation}
\begin{split}
[\vec{r}_1-\vec{r}, \vec{r}_2-\vec{r},\vec{r}_{\bar{1}\bar{2}}-\vec{r}]
&=[\vec{r}_1-\vec{r}, \vec{r}_2-\vec{r},\vec{r}_{\bar{1}\bar{2}}-\vec{r}_{\bar{2}}+\vec{r}_{\bar{2}}-\vec{r}]\\
&=[\vec{r}_1-\vec{r}, \vec{r}_2-\vec{r},-\frac{1}{\alpha_{\bar{1}\bar{2}}}(\vec{r}_{\bar{2}1}-\vec{r}_{\bar{2}})]\\
&=[\vec{r}_1-\vec{r}, \vec{r}_2-\vec{r},-\frac{1}{\alpha_{\bar{1}\bar{2}}}(\vec{r}_{\bar{1}2}-\vec{r}+\vec{r}-\vec{r}_{\bar{2}})]\\
&=-\frac{1}{b_{\bar{1}}}(a_{\bar{2}}-1)\frac{-1}{\alpha_{\bar{1}\bar{2}}}[\vec{r},\vec{r}_1,\vec{r}_2],
\end{split}
\end{equation}
and from it we obtain
\begin{equation}
 \alpha_{\bar{1}\bar{2}}<0.
\end{equation}
Finally, from
\begin{equation}
\begin{split}
[\vec{r}_1-\vec{r}, \vec{r}_2-\vec{r},\vec{r}_{\bar{1}\bar{1}}-\vec{r}]
&=[\vec{r}_1-\vec{r}, \vec{r}_2-\vec{r},\vec{r}_{\bar{1}\bar{1}}-\vec{r}_{\bar{1}}+\vec{r}_{\bar{1}}-\vec{r}]\\
&=[\vec{r}_1-\vec{r}, \vec{r}_2-\vec{r},\frac{\beta_{\bar{1}\bar{1}}}{\alpha_{\bar{1}\bar{1}}}(\vec{r}_{\bar{1}2}-\vec{r}_{\bar{1}})]\\
&=\frac{\beta_{\bar{1}\bar{1}}}{\alpha_{\bar{1}\bar{1}}}[\vec{r}_1-\vec{r}, \vec{r}_2-\vec{r},\vec{r}_{\bar{1}2}-\vec{r}+\vec{r}-\vec{r}_{\bar{1}}]\\
&=-\frac{1}{c_{\bar{1}}}(a_{\bar{1}}-1)\frac{\beta_{\bar{1}\bar{1}}}{\alpha_{\bar{1}\bar{1}}}[\vec{r},\vec{r}_1,\vec{r}_2],
\end{split}
\end{equation}
we find
\begin{equation}
 \frac{\beta_{\bar{1}\bar{1}}}{\alpha_{\bar{1}\bar{1}}}\geq0.
\end{equation}
Then it is easy to conclude
\begin{prop}\label{pro-convex}
A discrete surface is convex at the point $\vec{r}(n_1,n_2)$  if and only if
\begin{equation}
\begin{split}
\delta\geq 0,\quad \beta\geq0,\quad c_{\bar{1}}(a_{\bar{1}}-1)(a-1)<0,\quad b_{\bar{2}}(a_{\bar{2}}-1)(a-1)<0,
\end{split}
\end{equation}
$$$$
$$ \frac{\delta_{\bar{2}\bar{2}}}{\gamma_{\bar{2}\bar{2}}}\geq0,\quad \frac{\beta_{\bar{1}\bar{1}}}{\alpha_{\bar{1}\bar{1}}}\geq0,\quad \gamma_{\bar{1}\bar{2}}<0,\quad \alpha_{\bar{1}\bar{2}}<0.$$
\end{prop}
It is immediate to see from this proposition
\begin{cor}If a discrete centroaffine indefinite surface is locally convex everywhere, then $\delta=\beta=0,\gamma<0,\alpha<0$.
\end{cor}
By the Eq. (\ref{com-con-c}) and Proposition \ref{pro-convex}, we obtain
\begin{cor}A discrete centroaffine indefinite surface with constant coefficients is locally convex everywhere if and only if
\begin{equation}\label{C-Convex}
  \alpha=c<0, \gamma=b<0, \beta=\delta=0, a=b+c-bc.
\end{equation}
\end{cor}
\begin{rem}The discrete affine sphere shown in the left of Figure \ref{figure-AS} is a closed locally convex everywhere.
\end{rem}
We give another example which is local convex discrete centroaffine indefinite surface with constant coeffients.
\par
\begin{figure}[hbtp]
            \centering
            \begin{tabular}{cc}
              \includegraphics[width=.25\textwidth]{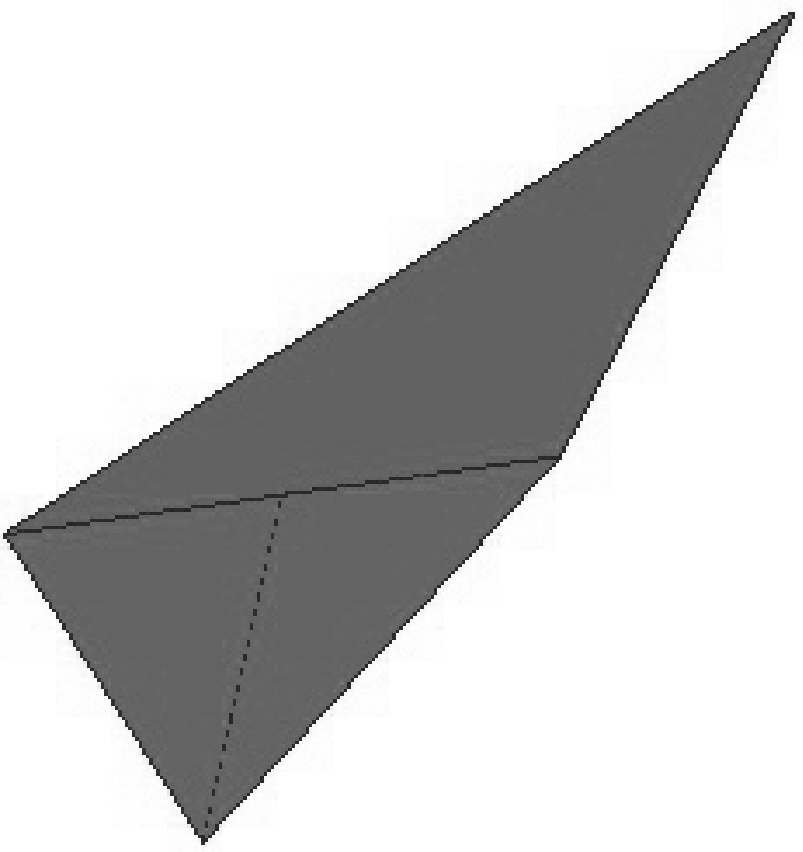} & \includegraphics[width=.25\textwidth]{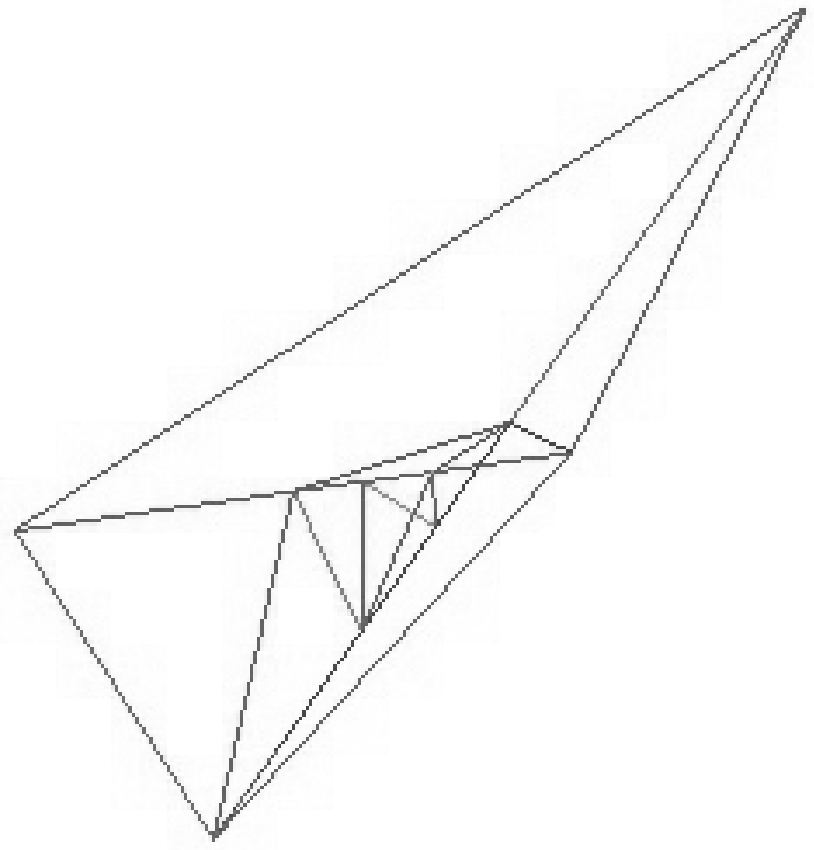}
            \end{tabular}
            \caption{Locally convex discrete centroaffine indefinite surface. }
            \label{figure-convex}
 \end{figure}

 {\bf Example.} It is convenient to check  $\alpha=c=-2, \gamma=b=-1, a=-5$ and $\beta=\delta=0$ satisfy Eq. (\ref{C-Convex}). Thus, we can obtain a discrete surface shown in Figure \ref{figure-convex} with 25 points. In fact, there are only 10 points because $\vec{r}_{22}=\vec{r}$.
\par If a discrete centroaffine indefinite surface is locally strongly convex, all determinants in Eq. (\ref{Det-all}) can not be zero. Hence we have
\begin{cor}There is not a discrete centroaffine indefinite surface which is locally strongly convex everywhere.
\end{cor}
{\bf Proof.} If a discrete surface is locally strongly convex everywhere, from the above calculations we have
$$\delta> 0,\quad \beta>0,\quad c_{\bar{1}}(a_{\bar{1}}-1)(a-1)<0,\quad b_{\bar{2}}(a_{\bar{2}}-1)(a-1)<0,$$
$$ \frac{\delta_{\bar{2}\bar{2}}}{\gamma_{\bar{2}\bar{2}}}>0,\quad \frac{\beta_{\bar{1}\bar{1}}}{\alpha_{\bar{1}\bar{1}}}>0,\quad \gamma_{\bar{1}\bar{2}}<0,\quad \alpha_{\bar{1}\bar{2}}<0,$$
which imply $\gamma_{\bar{2}\bar{2}}>0, \alpha_{\bar{1}\bar{1}}>0.$ This is inconsistent with $\gamma_{\bar{1}\bar{2}}<0,\quad \alpha_{\bar{1}\bar{2}}<0.$
\\ \rightline{$\Box$}

\renewcommand{\refname}{\bf\fontsize{12}{12}\selectfont References}
\bibliographystyle{amsplain}

\end{document}